\documentclass[11pt]{amsart}
\usepackage[margin=1in]{geometry}
\usepackage{graphicx} 
\usepackage{amsmath}
\usepackage{caption}
\usepackage{subfigure}
\usepackage{amsfonts}
\usepackage{dsfont}
\usepackage{amsthm}
\usepackage{geometry}
\usepackage{pgfplots}
\pgfplotsset{compat=1.18} 
\usepackage{color}
\usepackage{amsthm}
\usepackage{comment}
\usepackage[utf8]{inputenc}
\usepackage{hyperref}
\usepackage{bbm}
\usepackage{amssymb}

\newtheorem{theorem}{Theorem}[section]
\newtheorem{definition}{Definition}[section]
\newtheorem{lemma}{Lemma}[section]
\newtheorem{proposition}{Proposition}[section]

\newcommand{\Var}{\operatorname{Var}}

\newcommand{\1}{\mathds{1}}
\newcommand{\Prob}{\mathbb{P}}
\newcommand{\Ex}{\mathbb{E}}
\title[Descents and Inversions of Shelf-shuffles]{Limit theorems for descents and inversions of shelf-shuffles}

\author{Alexander Clay}
\date{\today}
\address{Department of Mathematics, University of Southern California (USC), Los Angeles, CA, USA.}
\email{ajclay@usc.edu}

\begin{document}
\begin{abstract}
   We prove central limit theorems for the number of descents and inversions of permutations produced by shelf-shuffles. These are a model for casino card shuffling machines. We show the asymptotic normality of the number of descents in two limiting regimes depending on the ratio of cards to shelves. On the other hand, we study the inversions by employing a modification of the techniques from Islak's analysis of the statistics of riffle shuffles. In particular, we obtain a bound for the rate of convergence for inversions that is independent of the number of shelves.
\end{abstract}
\maketitle

\section{Introduction and Main Results}
\label{Section1}
The mathematics of card shuffling is a rich subject that combines many aspects of probability, combinatorics, and algebra. For a good survey, we recommend Diaconis and Fulman \cite{MathofShufflingCards}. Most models of card shuffling analyzed in the literature are based on how humans shuffle cards. However, a certain mechanical card shuffling machine, known as a shelf-shuffler, has attracted significant mathematical interest. The shelf-shuffle was introduced by Diaconis, Fulman, and Holmes \cite{FDH} as a model for shelf-shuffling machines. 
\begin{definition}[Shelf-shuffle]
\label{def:shelfshuffle}
    Let $m\geq 1$ and $n\geq1$ be integers. Take $2m$ urns $U_1,\ldots,U_{2m}$. Place each of $n$ cards (labeled from $1$ to $n$) in an urn, independently of the other cards, with probability $1/(2m)$ in each urn. An $m$-shelf shuffle (of $n$ cards) is a random permutation $\sigma$ in the symmetric group $S_n$ obtained by reordering the cards in increasing (resp. decreasing) order for urns of odd (resp. even) index.
\end{definition}

For example, suppose $n=26$ and $m=2$, and we assign $26$ balls to $4$ urns as in Definition $\ref{def:shelfshuffle}$. Let $w$ be the random word as in \cite{Janson} with $w(i)$ equal to the index of the urn assigned to card $i$. In two-line notation, one possible sample would be
\begingroup
\setcounter{MaxMatrixCols}{40}   
\setlength{\arraycolsep}{1.6pt}  
\fontsize{9pt}{10pt}\selectfont                           
\begin{equation}
\label{eq:urnarrangement}
w=\begin{pmatrix}
  1 & 2 & 3 & 4 & 5 & 6 & 7 & 8 & 9 & 10 & 11 & 12 & 13 & 14 & 15 & 16 & 17 & 18 & 19 & 20 & 21 & 22 & 23 & 24 & 25 & 26 \\
  \color{blue}3 & \color{red}4 & \color{blue}1 & \color{blue}1 & \color{red}2 & \color{red}4 & \color{red}4 & \color{blue}3 & \color{blue}1 & \color{red}2 & \color{blue}3 & \color{blue}3 & \color{red}2 & \color{blue}1 & \color{red}4 & \color{red}2 & \color{blue}3 & \color{red}4 & \color{blue}1 & \color{red}2 & \color{red}4 & \color{blue}1 & \color{red}4 & \color{blue}1 & \color{blue}1 & \color{red}2 
\end{pmatrix}.
\end{equation}
\endgroup
If we sort the top row of $w$ according to the rules in Definition $\ref{def:shelfshuffle}$, we get
\begingroup
\setcounter{MaxMatrixCols}{40}   
\setlength{\arraycolsep}{1.6pt}  
\fontsize{9pt}{10pt}\selectfont                           
$$ 
\begin{pmatrix}
  3 & 4 & 9 & 14 & 19 & 22 & 24 & 25 & 26 & 20 & 16 & 13 & 10 & 5 & 1 & 8 & 11 & 12 & 17 & 23 & 21 & 18 & 15 & 7 & 6 & 2 \\
  \color{blue}1 & \color{blue}1 & \color{blue}1 & \color{blue}1 & \color{blue}1 & \color{blue}1 & \color{blue}1 & \color{blue}1 & \color{red}2 & \color{red}2 & \color{red}2 & \color{red}2 & \color{red}2 & \color{red}2 & \color{blue}{3} & \color{blue}3 & \color{blue}3 & \color{blue}3 & \color{blue}3 & \color{red}4 & \color{red}4 & \color{red}4 & \color{red}4 & \color{red}4 & \color{red}4 & \color{red}4
\end{pmatrix},
$$
\endgroup
and the resulting $2$-shelf-shuffle permutation is
\begingroup
\setcounter{MaxMatrixCols}{40}   
\setlength{\arraycolsep}{1.6pt}  
\fontsize{9pt}{10pt}\selectfont                           
\begin{equation} 
\label{eq:shelfshuffleintrodiag}
\sigma=\begin{pmatrix}
  1 & 2 & 3 & 4 & 5 & 6 & 7 & 8 & 9 & 10 & 11 & 12 & 13 & 14 & 15 & 16 & 17 & 18 & 19 & 20 & 21 & 22 & 23 & 24 & 25 & 26 \\
  \color{blue}3 & \color{blue}4 & \color{blue}9 & \color{blue}14 & \color{blue}19 & \color{blue}22 & \color{blue}24 & \color{blue}25 & \color{red}26 & \color{red}20 & \color{red}16 & \color{red}13 & \color{red}10 & \color{red}5 & \color{blue}1 & \color{blue}8 & \color{blue}11 & \color{blue}12 & \color{blue}17 & \color{red}23 & \color{red}21 & \color{red}18 & \color{red}15 & \color{red}7 & \color{red}6 & \color{red}2
\end{pmatrix}
\end{equation}
\endgroup
which is diagramatically represented in Figure $\ref{fig:shelfshufflediag}$.
\begin{figure}[htbp]
\centering
\begin{tikzpicture}[scale=0.2]
\def\n{26}
\draw (0,0) grid (\n,\n);
\fill[blue] (0,2) rectangle (1,3);
\fill[blue] (1,3) rectangle (2,4);
\fill[blue] (2,8) rectangle (3,9);
\fill[blue] (3,13) rectangle (4,14);
\fill[blue] (4,18) rectangle (5,19);
\fill[blue] (5,21) rectangle (6,22);
\fill[blue] (6,23) rectangle (7,24);
\fill[blue] (7,24) rectangle (8,25);
\fill[red] (8,25) rectangle (9,26);
\fill[red] (9,19) rectangle (10,20);
\fill[red] (10,15) rectangle (11,16);
\fill[red] (11,12) rectangle (12,13);
\fill[red] (12,9) rectangle (13,10);
\fill[red] (13,4) rectangle (14,5);
\fill[blue] (14,0) rectangle (15,1);
\fill[blue] (15,7) rectangle (16,8);
\fill[blue] (16,10) rectangle (17,11);
\fill[blue] (17,11) rectangle (18,12);
\fill[blue] (18,16) rectangle (19,17);
\fill[red] (19,22) rectangle (20,23);
\fill[red] (20,20) rectangle (21,21);
\fill[red] (21,17) rectangle (22,18);
\fill[red] (22,14) rectangle (23,15);
\fill[red] (23,6) rectangle (24,7);
\fill[red] (24,5) rectangle (25,6);
\fill[red] (25,1) rectangle (26,2);
\end{tikzpicture}
\caption{\label{fig:shelfshufflediag} Diagrammatic representation of the permutation $\sigma$ in Equation $\eqref{eq:shelfshuffleintrodiag}$. Each filled in box is at row $\sigma(i)$ and column $i$. \textcolor{blue}{Blue} (resp. \textcolor{red}{red}) boxes represent cards placed in \textcolor{blue}{odd} (resp. \textcolor{red}{even}) -indexed urns.}
\end{figure}

The authors in \cite{FDH} used the description from Definition $\ref{def:shelfshuffle}$ to construct the corresponding probability measure $\mathbb{P}_m$ on the symmetric group $S_n$ which gives the chance of each permutation of cards $\sigma$ after an $m$-shelf-shuffle. Additionally, they showed that $\mathbb{P}_m$ is close to the uniform distribution on $S_n$ in separation and $\ell_{\infty}$ distances when $m=O(n^{3/2})$. Letting $\Prob^{\ast k}_{m}$ be the induced measure on $S_n$ after $k$ iterated $m$-shelf shuffles, Chen and Ottolini \cite{ChenOtt} showed that there is cutoff in the total variation distance between $\Prob^{\ast k}_{m}$ and the uniform distribution when $k=\frac{5\log_{2m} n}{4}+O(1)$. Towards a different direction, the author has considered guessing strategies for shelf-shuffles with different levels of feedback \cite{Clay2025}, and his results were extended recently by Kuba \cite{kuba2026cardguessingsingleshelf} and Tripathi \cite{tripathi2026positionmatrixsingleshelfshuffle}. From a statistical standpoint, a limit theorem for the joint cycle structure of shelf-shuffles was proven in \cite{FDH} and convergence rates were found by Ren \cite{RenCycleStructure}. 

Descents and inversions are two more useful statistical tests of randomness in permutations. A permutation $\sigma\in S_n$ is said to have a descent at index $i$ if $\sigma(i)>\sigma(i+1)$, while an inversion is a pair $(i,j)$ with $i<j$ and $\sigma(i)>\sigma(j)$. Given a permutation $\sigma\in S_n$, we let $d(\sigma)$ and $I(\sigma)$ be the number of descents and the number of inversions, respectively. For example, the permutation 
\[\sigma=\begin{pmatrix} 1 & 2 & 3 & 4 & 5\\ 2&3&5&4&1\end{pmatrix}\] 
has $d(\sigma)=2$ and $I(\sigma)=5$. Diaconis and Fulman \cite{MathofShufflingCards} survey many interesting combinatorial and probabilistic results on the descents and inversions of riffle shuffles. 

In this paper, we study the descents and inversions of shelf-shuffles. Our first contribution is a central limit theorem for the number of descents in two asymptotic limiting regimes depending on the ratio of cards to shelves.
\begin{theorem}
\label{CLT descents}
Let $d(\sigma)$ be the number of descents of a permutation $\sigma$ obtained by an $m$-shelf shuffle of $n$ cards. We have the following.
\begin{enumerate}
    \item
        Suppose that $m\ll n^{1/2}$. Then, 
        \[\frac{d(\sigma)-n/2}{\sqrt{n/4}}\Rightarrow_d\mathcal{N}(0,1).\]
    \item Suppose that $m\gg n^{3/2}$. Then, 
    \begin{equation}
    \label{eq:cltdescentsn32}
        \frac{d(\sigma)-n/2}{\sqrt{n/12}}\Rightarrow_d\mathcal{N}(0,1).
    \end{equation}
\end{enumerate}
Here, $\Rightarrow_d$ denotes convergence in distribution.
\end{theorem}
 Before stating our result on inversions, we define the Kolmogorov distance between two random variables $X$ and $Y$ by $d_K(X,Y)=\sup_{t\in\mathbb{R}}|F_X(t)-F_Y(t)|$, where $F_X$ and $F_Y$ are the cumulative distribution functions. It is immediate that the Kolmogorov distance defines a metric on distributions, and that convergence in Kolmogorov distance implies convergence in distribution. Our second main result gives a central limit theorem with a convergence rate for the number of inversions.
\begin{theorem}
\label{CLT inversions}
    Let $I(\sigma)$ be the number of inversions of a permutation $\sigma$ obtained by an $m$-shelf shuffle of $n$ cards. Denoting $Z\sim\mathcal{N}(0,1)$, and letting 
    \begin{equation}
    \label{eq:expforynm}
        Y_{n,m}=\frac{\sqrt{n}}{2\binom{n}{2}\sqrt{\frac{m^2+2}{36m^2}}}\left(I(\sigma)-\frac{\binom{n}{2}}{2}\right),
    \end{equation}
    there exists an explicit constant $C$, independent of $n$ and $m$, such that $d_K(Y_{n,m},Z)\leq C/\sqrt{n}$.
\end{theorem}
The proof of Theorem $\ref{CLT inversions}$ uses methods from Islak \cite{Islak}, who proved central limit theorems for the number of descents and inversions of riffle shuffles. A $2m$-riffle shuffle utilizes the urn setup as in Definition $\ref{def:shelfshuffle}$, but the cards are ordered in increasing order in every urn, and then the inverse of the resulting permutation is taken. For example, if $w$ is the random word from Equation $\eqref{eq:urnarrangement}$, then sorting $w$ according to the rules for riffle shuffles gives
\begingroup
\setcounter{MaxMatrixCols}{40}   
\setlength{\arraycolsep}{1.6pt}  
\fontsize{9pt}{10pt}\selectfont  
$$ 
\pi=\begin{pmatrix}
  1 & 2 & 3 & 4 & 5 & 6 & 7 & 8 & 9 & 10 & 11 & 12 & 13 & 14 & 15 & 16 & 17 & 18 & 19 & 20 & 21 & 22 & 23 & 24 & 25 & 26 \\
  \color{blue}3 & \color{blue}4 & \color{blue}9 & \color{blue}14 & \color{blue}19 & \color{blue}22 & \color{blue}24 & \color{blue}25 & \color{red}5 & \color{red}10 & \color{red}13 & \color{red}16 & \color{red}20 & \color{red}26 & \color{blue}1 & \color{blue}8 & \color{blue}11 & \color{blue}12 & \color{blue}17 & \color{red}2 & \color{red}6 & \color{red}7 & \color{red}15 & \color{red}18 & \color{red}21 & \color{red}23
\end{pmatrix},
$$
\endgroup
and the resulting $4$-riffle shuffle permutation is
\begingroup
\setcounter{MaxMatrixCols}{40}   
\setlength{\arraycolsep}{1.6pt}  
\fontsize{9pt}{10pt}\selectfont                           
$$ 
\pi^{-1} = \begin{pmatrix}
  1 & 2 & 3 & 4 & 5 & 6 & 7 & 8 & 9 & 10 & 11 & 12 & 13 & 14 & 15 & 16 & 17 & 18 & 19 & 20 & 21 & 22 & 23 & 24 & 25 & 26 \\
  15 & 20 & 1 & 2 & 9 & 21 & 22 & 16 & 3 & 10 & 17 & 18 & 11 & 4 & 23 & 12 & 19 & 24 & 5 & 13 & 25 & 6 & 26 & 7 & 8 & 14
\end{pmatrix}.
$$
\endgroup
A permutation and its inverse have the same number of inversions, so we can modify Islak's methods to prove Theorem $\ref{CLT inversions}$. However, the number of descents is not preserved under taking inverses in general, so we developed other techniques to analyze the descents. With $\sigma$ a shelf-shuffle permutation, the authors in \cite{FDH} found that the mean and variance of the number of descents of $\sigma^{-1}$ are $(n-1)/2$ and $((n+1)/12)+(n-2)/(6m^2)$, respectively.

To connect our results with other literature on the mixing properties of the shelf-shuffle, it is useful to examine the relationship between descents, inversions, and valleys. A permutation $\sigma$ is said to have a valley at position $i$ if $\sigma(i-1)>\sigma(i)$ and $\sigma(i)<\sigma(i+1)$. In \cite{FDH}, it is shown that the number of valleys is a sufficient statistic to determine the chance of a shelf-shuffle permutation. Asymptotics for the number of valleys were used to prove the bounds on separation, $\ell_{\infty}$, and total variation distances in \cite{FDH} and \cite{ChenOtt}. Unfortunately, the descents and inversions of shelf-shuffles cannot be used to evaluate mixing in the same way. For example, the $2$-shelf shuffle permutation from Equation $\eqref{eq:shelfshuffleintrodiag}$ has $12$ descents and $166$ inversions, and the $10$-shelf shuffle permutation from Equation $\eqref{eq:sparseshelves}$ has $13$ descents and $166$ inversions. However, their global structures as shown in Figure $\ref{fig:shelfshufflediag}$ and Figure $\ref{fig:shelfshufflediag2}$, respectively, are dissimilar because of the difference in shelves. 

To explain the positive consequences of our results, we note that Theorem $\ref{CLT descents}$ shows that when $m\gg n^{3/2}$, the number of descents of a shelf-shuffle has the same limiting distribution as the descents of a uniformly random permutation. On the other hand, in Equation $\eqref{eq:expforynm}$, we see that $\Var(Y_{n,m})/n^3\sim(1/36)+(2/m^2)\to1/36$ when $m\to\infty$. This matches the asymptotic variance (divided by $n^3$) of the number of inversions of a uniformly random permutation. Therefore, the weak limit of the number of inversions in a shelf shuffle coincides with that for a uniformly random permutation as $n\to\infty$ provided that $m\to\infty$ as well.
\section*{Acknowledgments}
The author would like to thank his advisor, Jason Fulman, for suggesting this project and for his references. We also thank Evgeni Dimitrov and Raghavendra Tripathi for their kind advice.
\section{Descents}
\label{Section2}
The goal of this section is to construct a coupling for the number of descents which will be sufficient to prove Theorem $\ref{CLT descents}$. To motivate the appearance of two limiting regimes in Theorem $\ref{CLT descents}$, consider an example of a shelf-shuffle when the ratio of cards to shelves is small. Suppose $m=10$ and $n=26$, and we sample the random word 
    \begingroup
    \setcounter{MaxMatrixCols}{40}   
    \setlength{\arraycolsep}{1.6pt}  
    \fontsize{9pt}{10pt}\selectfont    
    $$ 
        w=\begin{pmatrix}
          1 & 2 & 3 & 4 & 5 & 6 & 7 & 8 & 9 & 10 & 11 & 12 & 13 & 14 & 15 & 16 & 17 & 18 & 19 & 20 & 21 & 22 & 23 & 24 & 25 & 26 \\
          \color{blue}7 & \color{red}10 & \color{red}18 & \color{blue}15 & \color{blue}19 & \color{red}10 & \color{red}8 & \color{red}6 & \color{blue}11 & \color{blue}19 & \color{red}2 & \color{blue}13 & \color{red}10 & \color{blue}5 & \color{blue}9 & \color{red}20 & \color{red}14 & \color{red}14 & \color{blue}7 & \color{red}12 & \color{red}4 & \color{blue}17 & \color{blue}3 & \color{blue}17 & \color{red}6 & \color{blue}19
        \end{pmatrix}.
    $$
    \endgroup
The resulting $10$-shelf shuffle permutation is 
\begingroup
\setcounter{MaxMatrixCols}{40}
\setlength{\arraycolsep}{1.6pt}
\fontsize{9pt}{10pt}\selectfont
\begin{equation}
\label{eq:sparseshelves}
\sigma=
\begin{pmatrix}
1 & 2 & 3 & 4 & 5 & 6 & 7 & 8 & 9 & 10 & 11 & 12 & 13 & 14 & 15 & 16 & 17 & 18 & 19 & 20 & 21 & 22 & 23 & 24 & 25 & 26 \\
\color{blue}{11} & \color{red}{23} & \color{blue}{21} & \color{red}{14} & \color{blue}{25} & \color{red}{8} & \color{blue}{1} & \color{blue}{19} & \color{red}{7} & \color{blue}{15} & \color{red}{13} & \color{red}{6} & \color{red}{2} & \color{blue}{9} & \color{red}{20} & \color{blue}{12} & \color{red}{18} & \color{red}{17} & \color{blue}{4} & \color{blue}{22} & \color{blue}{24} & \color{blue}{3} & \color{blue}{5} & \color{blue}{10} & \color{blue}{26} & \color{red}{16}
\end{pmatrix}.
\end{equation}
\endgroup
\begin{figure}[htbp]
\centering
\begin{tikzpicture}[scale=0.2]
\def\n{26}
\draw (0,0) grid (\n,\n);

\fill[blue] (0,10) rectangle (1,11);
\fill[red]  (1,22) rectangle (2,23);
\fill[blue] (2,20) rectangle (3,21);
\fill[red]  (3,13) rectangle (4,14);
\fill[blue] (4,24) rectangle (5,25);
\fill[red]  (5,7) rectangle (6,8);
\fill[blue] (6,0) rectangle (7,1);
\fill[blue] (7,18) rectangle (8,19);
\fill[red]  (8,6) rectangle (9,7);
\fill[blue] (9,14) rectangle (10,15);
\fill[red]  (10,12) rectangle (11,13);
\fill[red]  (11,5) rectangle (12,6);
\fill[red]  (12,1) rectangle (13,2);
\fill[blue] (13,8) rectangle (14,9);
\fill[red]  (14,19) rectangle (15,20);
\fill[blue] (15,11) rectangle (16,12);
\fill[red]  (16,17) rectangle (17,18);
\fill[red]  (17,16) rectangle (18,17);
\fill[blue] (18,3) rectangle (19,4);
\fill[blue] (19,21) rectangle (20,22);
\fill[blue] (20,23) rectangle (21,24);
\fill[blue] (21,2) rectangle (22,3);
\fill[blue] (22,4) rectangle (23,5);
\fill[blue] (23,9) rectangle (24,10);
\fill[blue] (24,25) rectangle (25,26);
\fill[red]  (25,15) rectangle (26,16);

\end{tikzpicture}
\caption{\label{fig:shelfshufflediag2} Diagrammatic representation of the permutation $\sigma$ in Equation $\eqref{eq:sparseshelves}$.}
\end{figure}

Note that a descent at $(i,i+1)$ of $\sigma$ occurs when the filled box in column $i$ is higher than the filled box in column $i+1$. Comparing Figure $\ref{fig:shelfshufflediag2}$ with Figure $\ref{fig:shelfshufflediag}$, where the ratio of cards to shelves is high, we see fundamental differences in the descent structure. In particular, most of the descents in Figure $\ref{fig:shelfshufflediag}$ occur with strings of red blocks, which represent cards placed in the same even-indexed urn. On the other hand, in Figure $\ref{fig:shelfshufflediag2}$, many descents occur between cards from different urns. This motivates the appearance of two limiting regimes in Theorem $\ref{CLT descents}$, which we now investigate.

The next proposition relates the descents of a shelf-shuffle to its underlying structure. Let $U(\sigma(i))$ be the index of the urn assigned to card $\sigma(i)$. We say that a pair $(i,i+1)$ is a \textit{concatenation} if $U(\sigma(i))\neq U(\sigma(i+1))$. For example, in the permutation displayed in Figure $\ref{fig:shelfshufflediag2}$, $(1,2)$ and $(21,22)$ are concatenations, while $(11,12)$ is not. We can classify descents based on the corresponding urn arrangement. Denote the integers $\{1,2,\ldots,m\}$ by $[m]$.
\begin{proposition}
\label{prop:concat_even_descents}
    If the pair $(i,i+1)$ is a descent of $\sigma$, then either
    \begin{enumerate}
        \item $(i,i+1)$ is a concatenation, or
        \item $U(\sigma(i))=U(\sigma(i+1))=2k$ for some $k\in[m]$.
    \end{enumerate}
    If condition (2) holds, then $(i,i+1)$ is a descent of $\sigma$.
    \begin{proof}
        Suppose for sake of contradiction that $\sigma(i)>\sigma(i+1)$ and $U(\sigma(i))=U(\sigma(i+1))=2k-1$ for some $k\in[m]$. From Definition $\ref{def:shelfshuffle}$, we know that cards assigned to odd-labeled urns are placed in increasing order, and therefore $\sigma(i)<\sigma(i+1)$, a contradiction. Now assume that condition (2) holds. Then cards $\sigma(i)$ and $\sigma(i+1)$ are placed in the same even-indexed urn, so by Definition $\ref{def:shelfshuffle}$ they are placed next to each other in descending order. We conclude that $(i,i+1)$ is a descent of $\sigma$.
    \end{proof}
\end{proposition}
We say that a pair $(i,i+1)$ is a \textit{concatenated descent} (resp. \textit{even descent}) if it is a descent and condition (1) (resp. (2)) of Proposition $\ref{prop:concat_even_descents}$ holds. In Figure $\ref{fig:shelfshufflediag2}$, we have that $(2,3)$ is a concatenated descent, and $(11,12)$ is an even descent. Let $f(\sigma)$ (resp. $g(\sigma)$) be the number of concatenated (resp. even) descents, and let $d(\sigma)$ be the number of descents. Proposition $\ref{prop:concat_even_descents}$ implies that $d(\sigma)=f(\sigma)+g(\sigma)$. 

To get a better understanding of the random variables $f(\sigma)$ and $g(\sigma)$, we introduce more concepts. We say that an index $i$ is \textit{lonely} if $U(\sigma(i))\neq U(\sigma(j))$ for all $j\neq i$, i.e. card $\sigma(i)$ is the only card in its urn. In Figure $\ref{fig:shelfshufflediag2}$, we have that $1,2$ and $19$ are examples of lonely indices, while $11$ and $24$ are not.

It turns out that we can use lonely indices to couple shelf-shuffles with uniformly random permutations. For a subset $A\subseteq S_n$, let $S_A=\{\pi\in S_n:\pi(a)=a\;\;\text{for all } a\in A\}$ be the set of permutations $\pi\in S_n$ which fix all elements of $A$. The following lemma establishes the coupling which will be used to prove Theorem $\ref{CLT descents}.2$.
\begin{lemma}
\label{lem:lonelycoupling}
    Let $\sigma$ be a shelf-shuffle permutation, and let $L$ be its set of lonely indices. For any subset $A\subseteq L$, if we choose $\pi\in S_A$ uniformly at random and independently of $\sigma$, then $\sigma\circ\pi$ is a uniformly random permutation in $S_n$.
    \begin{proof}
        Let $A\subseteq L$, and let $k=|A|$. For an injection $\varphi:A\to[n]$, define the event
        \[E_{\varphi}=\{\sigma(a)=\varphi(a)\text{ for all } a\in A\}.\]
        We will show that all events $E_{\varphi}$ have the same probability regardless of the choice of injection $\varphi$. Let $\varphi,\psi:A\to[n]$ be injections. Consider a realization of the urn assignments producing a shelf-shuffle $\sigma$ in $E_{\varphi}$. For each $a\in A$, the card $\varphi(a)=\sigma(a)$ is lonely since $A\subseteq L$. Hence, $\varphi(a)$ is the unique card in its urn. 
        
        Replace the label $\varphi(a)$ by $\psi(a)$ for every $a\in A$, leaving all urn assignments unchanged and leaving all other card labels unchanged. Since each $\varphi(a)$ is the only card in its urn, the relabeling does not affect the shelf-shuffle ordering; the card occupying position $a$ in the shuffled deck is simply relabeled from $\varphi(a)$ to $\psi(a)$. The resulting shelf-shuffle permutation $\sigma'$ hence satisfies $\sigma'(a)=\psi(a)$ for all $a\in A$. This is a reversible construction by replacing $\psi(a)$ with $\varphi(a)$, so it gives a bijection between the outcomes in $E_{\varphi}$ and $E_{\psi}$. Since all urn assignment words are equally weighted, this implies that $\mathbb{P}(E_{\varphi})=\mathbb{P}(E_{\psi})$. The events $E_{\varphi}$ partition the sample space, and since there are $n!/(n-k)!:=(n)_k$ injective maps $A\to[n]$, we have $\mathbb{P}(E_{\varphi})=1/(n)_k$ for every injection $\varphi$.
        
        Let $\pi$ be uniformly distributed on $S_A$ and independent of $\sigma$. Fix $\tau\in S_n$. Then, 
        \begin{equation}
            \label{eq:condlexp}    \Prob(\sigma\circ\pi=\tau)=\sum_{\varphi}\Prob(\sigma\circ\pi=\tau\mid E_{\varphi})\Prob(E_{\varphi})
        \end{equation}
        where the sum runs over all injections $\varphi:A\to[n]$. If $E_{\varphi}$ happens, then $(\sigma\circ\pi)|_A=\varphi$ since each $\pi\in S_A$ fixes $A$ pointwise. Therefore, $\Prob(\sigma\circ\pi=\tau\mid E_{\varphi})=0$ unless $\varphi=\tau|_A$. If $\varphi=\tau|_A$, the map $\pi\to\sigma\circ\pi$ is a bijection from $S_A$ onto the set of permutations whose restriction to $A$ equals $\varphi$. Therefore, conditional on $E_{\varphi}$, the random permutation $\sigma\circ\pi$ is uniform among the $(n-k)!$ permutations extending $\varphi$. Hence,
        \begin{equation}
        \label{eq:etaua}
            \Prob(\sigma\circ\pi=\tau\mid E_{\tau\mid A})=\frac{1}{(n-k)!}
        \end{equation}
        Combining Equations $\eqref{eq:condlexp}$ and $\eqref{eq:etaua}$, we obtain that 
        \[\Prob(\sigma\circ\pi=\tau)=\frac{1}{(n)_k}\frac{1}{(n-k)!}=\frac{1}{n!}.\]
        Since this holds for every $\tau\in S_n$, the lemma follows.
    \end{proof}
\end{lemma}
The next proposition makes the qualitative differences between the descent structure in Figure $\ref{fig:shelfshufflediag}$ and Figure $\ref{fig:shelfshufflediag2}$ precise. In particular, we show that when the ratio of cards to shelves is small (resp. large), most descents are concatenated (resp. even).
\begin{proposition}
\label{prop:descentsfg}
Let $\sigma$ be an $m$-shelf-shuffle of $n$ cards, and write the number of descents $d(\sigma)=f(\sigma)+g(\sigma)$, where $f(\sigma)$ (resp. $g(\sigma)$) is the number of concatenated (resp. even) descents. The random variables $f(\sigma)$ and $g(\sigma)$ satisfy the following properties.
\begin{enumerate}
    \item $f(\sigma)\leq 2m-1$ almost surely (a.s.).
    \item When $m\gg n^{3/2}$, we have that $\frac{f(\sigma)-n/2}{\sqrt{n/12}}\Rightarrow_d\mathcal{N}(0,1)$ as $n\to\infty$.
    \item If $m\gg n^{3/2}$, we have that $n^{-1/2}\cdot g(\sigma)\to 0$ in probability. 
    \item The convergence $\frac{g(\sigma)-n/2}{\sqrt{n/4}}\Rightarrow_d\mathcal{N}(0,1)$ as $n\to\infty$ holds when $m\ll n^{1/2}$.
\end{enumerate}
\end{proposition}
Before proving the proposition, we develop some more theory. We call a pair of indices $(i,i+1)$ \textit{special} if both $i$ and $i+1$ are lonely. In Figure $\ref{fig:shelfshufflediag2}$, the pair $(2,3)$ is special. Let $f_{SP}(\sigma)$ be the number of special pairs of indices. Call a descent $(i,i+1)$ a $\textit{special concatenated descent}$ (SCD) if $(i,i+1)$ is special. We note that if $(i,i+1)$ is special and a descent then it is also an SCD. Letting $f_{SCD}(\sigma)$ be the number of SCDs, we have the following useful estimates, which will imply Proposition $\ref{prop:descentsfg}$.
\begin{lemma}
\label{lem:estimatesfordescents}
    Suppose $\sigma$ is an $m$-shelf-shuffle of $n$ cards. Then, the following hold almost surely.
    \begin{enumerate}
        \item $0\leq f(\sigma)-f_{SCD}(\sigma)\leq n-1-f_{SP}(\sigma)$.
        \item Let $A$ be the set of indices $i$ such that $(i-1,i)$ or $(i,i+1)$ is special (possibly both). For any $\pi\in S_A$, we have $0\leq d(\sigma\circ\pi)-f_{SCD}(\sigma)\leq n-1-f_{SP}(\sigma)$. 
    \end{enumerate}
    \begin{proof}
        We first prove $(1)$. Notice that $f(\sigma)-f_{SCD}(\sigma)$ is the number of pairs $(i,i+1)$ which are both concatenated descents and not special. This combinatorial description implies that $f(\sigma)-f_{SCD}(\sigma)$ is nonnegative and bounded above by the number of non-special pairs, proving $(1)$. To show $(2)$, first assume that $(i,i+1)$ is an SCD of $\sigma$; then $i\in A$ and $i+1\in A$ by special-ness, and so $(\sigma\circ\pi)(i)=\sigma(i)$ and $(\sigma\circ\pi)(i+1)=\sigma(i+1)$ because $\pi$ fixes elements of $A$. Thus, $(i,i+1)$ is a descent of $\sigma\circ\pi$, implying that $f_{SCD}(\sigma)\leq d(\sigma\circ\pi)$. On the other hand, suppose that $(i,i+1)$ is a descent of $\sigma\circ\pi$ which is not an SCD of $\sigma$. Assume for sake of contradiction that $(i,i+1)$ is special. Since $\pi$ fixes elements of $A$, we have that $(i,i+1)$ is a descent of $\sigma$ and moreover, by specialness, an SCD, which is a contradiction. Hence, $(i,i+1)$ is not special, so we conclude $(2)$.
    \end{proof}
\end{lemma}
\begin{proof}[Proof of Proposition $\ref{prop:descentsfg}$]
We first prove (1). Let $\sigma$ be a shelf-shuffle. Since there are $2m$ urns, there are at most $2m-1$ concatenations, and therefore at most $2m-1$ concatenated descents, so $f(\sigma)\leq 2m-1$ a.s. 

To establish (2), we let $A$ be the set of indices $i$ such that $(i-1,i)$ or $(i,i+1)$ is special (possibly both). Choose $\pi\in S_A$ uniformly at random. Noting that $A\subseteq L$, where $L$ is the set of lonely indices of $\sigma$, Lemma $\ref{lem:lonelycoupling}$ implies that $\sigma\circ\pi$ is a uniformly random permutation. Since $\sigma\circ\pi$ is uniformly random, it is well-known (see, e.g. \cite{8a557c1d-9f79-3773-9c3e-72728b184f9b}) that the number of descents $d(\sigma\circ\pi)$ satisfies a central limit theorem:
    \begin{equation}
    \label{eq:uniformdescents}
        \frac{d(\sigma\circ\pi)-n/2}{\sqrt{n/12}}\Rightarrow_d\mathcal{N}(0,1).
    \end{equation}
We will now show that $n^{-1/2}\cdot|d(\sigma\circ\pi)-f(\sigma)|\to 0$ in $L^1$ when $m\gg n^{3/2}$. From the triangle inequality, monotonicity, and Lemma $\ref{lem:estimatesfordescents}$, we get
\begin{equation}
    \label{eq:triangleineqforspecialdescents}
        \begin{aligned}
        \mathbb{E}\left|d(\sigma\circ\pi)-f(\sigma)\right|&\leq\Ex\left|d(\sigma\circ\pi)-f_{SCD}(\sigma)\right|+\Ex\left|f(\sigma)-f_{SCD}(\sigma)\right|\\ &\leq 2n-2\cdot\mathbb{E}[f_{SP}(\sigma)].\\
        \end{aligned}
\end{equation}
Therefore, we need to compute the expected number of special pairs $f_{SP}(\sigma)$. Let $E_i$ be the event that $(i,i+1)$ is special. Note that $E_i$ occurs if and only if card $\sigma(i+1)$ avoids urn $U(\sigma(i))$ and all other $n-2$ cards (besides cards $\sigma(i)$ and $\sigma(i+1)$) avoid both urns $U(\sigma(i))$ and $U(\sigma(i+1))$. Therefore, by the method of indicators,
    \begin{equation}
    \label{eq:expectedvalueofspecialdescents}
        \mathbb{E}[f_{SP}(\sigma)]=\sum_{i=1}^{n-1}\mathbb{P}(E_i)=\sum_{i=1}^{n-1}(1-1/2m)(1-1/m)^{n-2}=(n-1)(1-1/2m)(1-1/m)^{n-2}.
    \end{equation}
    Combining Equation $\eqref{eq:expectedvalueofspecialdescents}$ with the bound from Equation $\eqref{eq:triangleineqforspecialdescents}$, we obtain that
    \begin{equation}
    \label{eq:descentsboundonexp}
            \mathbb{E}\left|d(\sigma\circ\pi)-f(\sigma)\right|\leq 2(n-1)\left(1-(1-1/2m)(1-1/m)^{n-2}\right).
    \end{equation}
    Assuming that $m,n\to\infty$ in such a way that $n/m\to 0$ (which is satisfied when $m\gg n^{3/2}$), a binomial expansion yields
    \begin{equation}
    \label{eq:binomial}
        \begin{aligned}
            (1-1/2m)(1-1/m)^{n-2}&=(1-1/2m)\left(1-\frac{n-2}{m}+O\left(\frac{n^2}{m^2}\right)\right)\\
            &=1-\frac{n-2}{m}-\frac{1}{2m}+\frac{n-2}{2m^2}+O\left(\frac{n^2}{m^2}\right)\\
            &=1-\frac{n}{m}+O\left(\frac{n^2}{m^2}\right).
        \end{aligned}
    \end{equation}
    Substituting the estimate from Equation $\eqref{eq:binomial}$ into Equation $\ref{eq:descentsboundonexp}$ and multiplying by $n^{-1/2}$ gives
    \begin{equation}
        \begin{aligned}
            n^{-1/2}\cdot\mathbb{E}\left|d(\sigma\circ\pi)-f(\sigma)\right|&\leq\frac{2(n-1)}{\sqrt{n}}\left(\frac{n}{m}+O\left(\frac{n^2}{m^2}\right)\right)\\
            &=\frac{2n^{3/2}}{m}+O\left(\frac{n^{5/2}}{m^2}\right)
        \end{aligned}
    \end{equation}
    which tends to zero as $n\to\infty$ when $m\gg n^{3/2}$. This shows that $n^{-1/2}\cdot|d(\sigma\circ\pi)-f(\sigma)|\to 0$ in $L^1$ when $m\gg n^{3/2}$. This convergence, together with Equation $\eqref{eq:uniformdescents}$ and an application of Slutsky's theorem as in \cite[Theorem 2.7]{AsymptoticStatistics}, implies part (2).

    Next, we prove part (3). Let $g_C(\sigma)$ be the number of cards placed into even-labeled urns, and let $g_U(\sigma)$ be the number of occupied even-labeled urns. An occupied even-labeled urn contributes $k$ descents if and only if it is occupied with $k+1$ cards. This is because, according to Definition $\ref{def:shelfshuffle}$, cards in even-labeled urns are placed in descending order. Therefore, we have the identity $g(\sigma)=g_C(\sigma)-g_U(\sigma)$. Note that since each card is assigned independently with probability $1/2$ of being placed into an even-labeled bin, we have that $g_C(\sigma)\sim\text{Binomial}(n,1/2)$. The classical central limit theorem implies that
    \begin{equation}
    \label{eq:cltforgc}
        \frac{g_C(\sigma)-n/2}{\sqrt{n/4}}\Rightarrow_d\mathcal{N}(0,1)
    \end{equation}
    as $n\to\infty$. Since there are exactly $m$ even-labeled bins, we have that at most $m$ can be occupied, so $g_U(\sigma)\leq m$ a.s. Therefore, $n^{-1/2}\cdot g_U(\sigma)\to 0$ in probability as $n\to\infty$ when $m\ll n^{1/2}$. We conclude part (3) from Slutsky's theorem and the convergence in Equation $\eqref{eq:cltforgc}$.
    
    The last step is to establish (4). To show that $n^{-1/2}\cdot g(\sigma)\to 0$ in probability when $m\gg n^{3/2}$, we will show that $n^{-1/2}\cdot g(\sigma)\to 0$ in $L^1$ in the same limiting regime. Noting that $g_C(\sigma)\geq g_U(\sigma)$ a.s., we get $\Ex[g(\sigma)]=\Ex[g_C(\sigma)]-\Ex[g_U(\sigma)]$. In the proof of part (3), we established that $g_C(\sigma)\sim\text{Binomial}(n,1/2)$, so $\Ex[g_C(\sigma)]=n/2$. Computing the expected number of occupied even-labeled bins $\Ex[g_U(\sigma)]$, we let $E_i$ be the event that bin $i$ is occupied; then by considering unoccupied bins we have $\mathbb{P}(E_i)=1-\mathbb{P}(E_i^c)=1-(1-1/2m)^n$. By the method of indicators, we obatin
    \begin{equation}
    \label{eq:expectationofgsigma}
        \begin{aligned}
            \mathbb{E}|g(\sigma)|&=\mathbb{E}[g_C(\sigma)]-\Ex[g_U(\sigma)]\\
            &=\frac{n}{2}-\sum_{i=1}^{m}\mathbb{P}(E_{2i})\\
            &=\frac{n}{2}-m+m(1-1/2m)^n.\\
        \end{aligned}
    \end{equation}
    Substituting a binomial expansion for $(1-1/2m)^{n}$ into Equation $\eqref{eq:expectationofgsigma}$ and multiplying by $n^{-1/2}$ yields
    \begin{equation}
    \label{eq:binomialgsigma}
        \begin{aligned}
          n^{-1/2}\cdot\mathbb{E}|g(\sigma)|&=n^{-1/2}\left(\frac{n}{2}-m+m\left(1-\frac{n}{2m}+\frac{\binom{n}{2}}{4m^2}+O\left(\frac{n^3}{m^3}\right)\right)\right)\\
          &=\frac{n^{3/2}}{8m}+O\left(\frac{n^{5/2}}{m^2}\right)
        \end{aligned}
    \end{equation}
    which tends to zero when $m\gg n^{3/2}$ and $n\to\infty$. As in Equation $\eqref{eq:binomial}$, the estimates in the big-O notation in Equation $\ref{eq:binomialgsigma}$ are valid when $n/m\to 0$, which is true when $m\gg n^{3/2}$. Since convergence in $L^1$ implies convergence in probability, part (4) follows.
    \end{proof}
    \begin{proof}[Proof of Theorem $\ref{CLT descents}$]
    From Proposition $\ref{prop:descentsfg}$, we have the decomposition into concatenated and even descents: $d(\sigma)=f(\sigma)+g(\sigma)$. Hence, we may rewrite
    \[\frac{d(\sigma)-n/2}{\sqrt{n/4}}=\frac{g(\sigma)-n/2}{\sqrt{n/4}}+\frac{f(\sigma)}{\sqrt{n/4}}.\]
    By Proposition $\ref{prop:descentsfg}.1$, we have $n^{-1/2}\cdot f(\sigma)\to 0$ in probability when $m\ll n^{1/2}$. Slutsky's theorem along with the limiting distribution of $g(\sigma)$ in Proposition $\ref{prop:descentsfg}.4$ yields Theorem $\ref{CLT descents}.1$. 

    Next, we can rewrite
    \[\frac{d(\sigma)-n/2}{\sqrt{n/12}}=\frac{f(\sigma)-n/2}{\sqrt{n/12}}+\frac{g(\sigma)}{\sqrt{n/12}}.\]
    By Proposition $\ref{prop:descentsfg}.3$, we have $n^{-1/2}\cdot g(\sigma)\to 0$ in probability when $m\gg n^{3/2}$. Combining Slutsky's theorem with the weak convergence of $f(\sigma)$ in Proposition $\ref{prop:descentsfg}.2$ gives Theorem $\ref{CLT descents}.2$.
    \end{proof}
\section{Inversions}
\label{Section3}
In this section, we prove Theorem $\ref{CLT inversions}$. First, we collect a result from Chen and Shao \cite{ChenShao}. To state it precisely, we define the notion of a $U$-statistic.
\begin{definition}[$U$-statistics] A random variable $W_n$ is an order $r$ $U$-statistic with symmetric kernel $h$ if we can write
\begin{equation}
\label{def:U-statistic}
    W_n=\frac{1}{\binom{n}{r}}\sum_{1\leq i_1<i_2<\ldots<i_r\leq n}h(Z_{i_1},Z_{i_2},\ldots,Z_{i_r})
\end{equation}
where $h$ is a real-valued Borel measurable symmetric function and $Z_1,Z_2,\ldots,Z_n$ are independent and identically distributed (i.i.d.) random variables.
\end{definition}
Before stating a convergence result on $U$-statistics, we need some notation. Given a $U$-statistic $W_n$ with representation as in Equation $\eqref{def:U-statistic}$, let us set $h_1(Z_1)=\mathbb{E}[h(Z_1,\ldots,Z_r)$$\mid Z_1]$, $\sigma_1^2=\text{Var}(h_1(Z_1))$, $\sigma^2=\text{Var}(h(Z_1,\ldots,Z_r))$, $\zeta=\mathbb{E}|h_1(Z_1)|^3$. Also, let $Z\sim\mathcal{N}(0,1)$. The following is a special case of a result of Chen and Shao \cite{ChenShao}.
\begin{proposition}[{\cite[Theorem 3.1]{ChenShao}}] 
\label{prop:U-statistics}
Let $W_n$ be an order $r$ $U$-statistic with symmetric kernel $h$ and decomposition as in Equation $\eqref{def:U-statistic}$. Suppose that $\mathbb{E}[h(Z_1,\ldots,Z_r)]=0$, $\sigma^2<\infty$, $\sigma_1^2>0$, and $\zeta<\infty$. Then, 
\[d_K\left(\frac{\sqrt{n}}{r\sigma_1}W_n,Z\right)\leq \frac{7\cdot\zeta}{\sqrt{n}\sigma_1^3}+\frac{(1+\sqrt{2})(r-1)\sigma}{(r(n-r+1))^{1/2}\sigma_1}.\]
\end{proposition}
Next, we follow methods of Islak \cite{Islak} and develop a $U$-statistic framework for the inversions of shelf-shuffles. Let $X_i$ be the label of the urn assigned to card $i$. Islak showed that if $\pi$ is a $2m$-riffle shuffle permutation, then the number of inversions $I(\pi)$ satisfies $I(\pi)=_d\sum_{i<j}\1(X_i>X_j)$, which allows for $I(\pi)$ to be represented as a $U$-statistic. Here, $\1(\cdot)$ denotes an indicator. It turns out that we have an analogous description for the inversions of shelf-shuffles.
\begin{proposition}
\label{prop:invshelfshuffledecomp}
    Let $\sigma$ be an $m$-shelf-shuffle of $n$ cards, and let $X_1,\ldots,X_n$ be i.i.d. random variables uniformly distributed on $[2m]$. Let $A_{i,j}$ be the event that $X_i=X_j=2k$ for some $k\in [m]$. The number of inversions $I(\sigma)$ satisfies
    \[I(\sigma)=_d\sum_{1\leq i<j\leq n}\left(\1(X_i>X_j)+\1(A_{i,j})\right).\]
    \begin{proof}
        For any permutation $\sigma\in S_n$, we have that $I(\sigma)=I(\sigma^{-1})$. In particular, we have
        \[I(\sigma)=\sum_{1\leq i<j\leq n}\1(\sigma(i)>\sigma(j))=\sum_{1\leq i<j\leq n}\1(\sigma^{-1}(i)>\sigma^{-1}(j)).\]
        Now let $\sigma$ be a shelf-shuffle permutation. Note that for each $i\in[n]$, $\sigma^{-1}(i)$ is the position of card $i$ in the shuffled deck. Fixing $i<j$, we see from the construction of shelf-shuffles in Definition $\ref{def:shelfshuffle}$ that the position of card $i$ is greater than the position of card $j$ (i.e. $\sigma^{-1}(i)>\sigma^{-1}(j)$) if and only if one of the following two mutually exclusive possibilities happens.
        \begin{enumerate}
            \item Card $i$ is placed in a urn with strictly greater index than the urn card $j$ is placed into. This follows because the cards are sorted in blocks relative to the indices of their urns.
            \item Card $i$ and card $j$ are placed into the same even-indexed urn. This is because cards in even-indexed urns are placed in decreasing order relative to their labels.
        \end{enumerate}
        Letting $X_i$ be the index of the urn card $i$ is placed into, we have that the random variables  $X_1,\ldots,X_n$ are independent and uniformly distributed on $[2m]$ because each urn is equally likely. The proposition follows after noticing that possibility (1) corresponds to the event $X_i>X_j$ and possibility (2) corresponds to the event $A_{i,j}$.
    \end{proof}
\end{proposition}
We show that the structure in Proposition $\ref{prop:invshelfshuffledecomp}$ is sufficient to write the inversions as a $U$-statistic. 
\begin{proposition}
\label{prop:inversionsmaindecomp}
    Let $\sigma$ be an $m$-shelf shuffle of $n$ cards, and let $I(\sigma)$ be the number of inversions. Let $X_1,\ldots, X_n$ be independent and uniformly distributed over $[2m]$, and let $U_1,U_2,\ldots,U_n$ be i.i.d. continuous uniform random variables on $(0,1)$. Let $Z_i=(X_i,U_i)$ for each $i\in[n]$, and let the event $A_{i,j}=\{X_i=X_j,\text{ (both) even}\}$. Then, setting 
    \begin{equation}
    \label{eq:gzizj}
        g(Z_i,Z_j)=\binom{n}{2}\bigg[\left(\1(X_{i}>X_{j})+\1(A_{i,j})\right)\1(U_{i}<U_{j})+\left(\1(X_{i}<X_{j})+\1(A_{j,i})\right)\1(U_{i}>U_{j})\bigg]
    \end{equation}
    we have the distributional equality
    \[I(\sigma)=_d\frac{1}{\binom{n}{2}}\sum_{1\leq i<j\leq n}g(Z_i,Z_j).\]
\end{proposition}
\begin{proof}
Let the random variables $U_1,U_2,\ldots,U_n$ and $X_1,X_2,\ldots,X_n$ be as in the statement of the proposition. Order $U_1,U_2,\ldots,U_n$ as $U_{\tau(1)}<U_{\tau(2)}<\cdots<U_{\tau(n)}$, where $\tau\in S_n$ is the appropriate permutation. Then $\tau$ is uniformly distributed over $S_n$, so we have
\begin{equation}
\label{eq:equalityofxi}
(X_1,\ldots,X_n)=_d(X_{\tau(1)},\ldots,X_{\tau(n)}).
\end{equation}
where $=_d$ denotes equality in distribution. It follows from Proposition $\ref{prop:invshelfshuffledecomp}$ and Equation $\eqref{eq:equalityofxi}$ that
\begin{equation*}
    \begin{aligned}
        I(\sigma)&=_d \sum_{i<j}\left(\1(X_{\tau(i)}>X_{\tau(j)})+\1(A_{\tau(i),\tau(j)})\right)\\
        &=_d\sum_{i,j=1}^n\left(\1(X_{\tau(i)}>X_{\tau(j)},i<j)+\1(A_{\tau(i),\tau(j)},i<j)\right).
    \end{aligned}
\end{equation*}
Noting that $i<j$ if and only if $U_{\tau(i)}<U_{\tau(j)}$, it follows that
\begin{equation}
    \begin{aligned}
    I(\sigma)&=_d\sum_{i,j=1}^n\left(\1(X_{\tau(i)}>X_{\tau(j)},U_{\tau(i)}<U_{\tau(j)})+\1(A_{\tau(i),\tau(j)},U_{\tau(i)}<U_{\tau(j)})\right)\\
    &=\sum_{i,j=1}^n\left(\1(X_{i}>X_{j})\1(U_i<U_j)+\1(A_{i,j})\1(U_{i}<U_{j})\right)\\
    &=\sum_{i,j=1}^n\left(\1(X_{i}>X_{j})+\1(A_{i,j})\right)\1(U_{i}<U_{j}).
    \end{aligned}
\end{equation}
Set $Z_i=(X_i,U_i)$, and note that the random variables $Z_1,Z_2,\ldots, Z_n$ are i.i.d. Define the function $A(x_i,x_j)=\1(x_i=x_j,\text{ (both) } x_i, x_j\text{ even})$, and note that $A$ is symmetric and $A(X_i,X_j)=A_{i,j}$. Then set
\[f((x_i,u_i),(x_j,u_j))=\binom{n}{2}\left(\1(x_{i}>x_{j})+A(x_i,x_j)\right)\1(u_{i}<u_{j})\]
and let 
\[g((x_i,u_i),(x_j,u_j))=f((x_i,u_i),(x_j,u_j))+f((x_j,u_j),(x_i,u_i)).\]
Note that $g(Z_i,Z_j)$ matches Equation $\eqref{eq:gzizj}$. By symmetry in the indices of summation, we have
\[I(\sigma)=_d\sum_{i,j=1}^n\left(\1(X_{i}>X_{j})+\1(A_{i,j})\right)\1(U_{i}<U_{j})=\frac{1}{\binom{n}{2}}\sum_{i<j}g(Z_i,Z_j)\]
as desired.
\end{proof}
Since Proposition $\ref{prop:U-statistics}$ only applies to mean zero $U$-statistics, we need to center $I(\sigma)$ by its mean. From Proposition $\ref{prop:invshelfshuffledecomp}$ and symmetry, we find that
\[\mathbb{E}[I(\sigma)]=\binom{n}{2}\left(\mathbb{P}(X_1>X_2)+\mathbb{P}(A_{1,2})\right).\]
By independence and symmetry, we have 
\[\mathbb{P}(X_1>X_2)=(1-\mathbb{P}(X_1=X_2))/2=\frac{1}{2}-\frac{1}{4m}.\]
Moreover, we have that
\[\mathbb{P}(A_{1,2})=\mathbb{P}(X_1=X_2,\text{ both even})=\mathbb{P}(X_1=X_2)\mathbb{P}(X_1\text{ even})=\frac{1}{4m}.\]
It follows that $\mathbb{E}[I(\sigma)]=\binom{n}{2}/2$. Letting 
\begin{equation}
\label{eq:hzizj}
        h(Z_i,Z_j)=g(Z_i,Z_j)-\frac{\binom{n}{2}}{2}
\end{equation}
we will collect the ingredients to apply Proposition $\ref{prop:U-statistics}$ to the random variable
\begin{equation}
\label{eq:wn}
    W_n:=\frac{1}{\binom{n}{2}}\sum_{1\leq i<j\leq n}h(Z_i,Z_j)=_dI(\sigma)-\frac{\binom{n}{2}}{2}.
\end{equation}
\begin{proposition}
\label{prop:dataforclt}
Let $g(Z_i,Z_j)$ be as in Equation $\eqref{eq:gzizj}$, and set $h(Z_i,Z_j)=g(Z_i,Z_j)-\binom{n}{2}/2$. Denote $h_1(Z_1)=\Ex[h(Z_1,Z_2)\mid Z_1]$. We have the following.
\begin{enumerate}
    \item $\Ex[h(Z_1,Z_2)]=0$.
    \item $\sigma^2:=\Var(h(Z_1,Z_2))=\binom{n}{2}^2/4$.
    \item $\sigma_1^2:=\Var(h_1(Z_1))=\frac{m^2+2}{36m^2}$.
    \item $\zeta:=\mathbb{E}|h_1(Z_1)|^3\leq\left(\frac{5}{2}\right)^3\binom{n}{2}^3$.
\end{enumerate}
\end{proposition}
\begin{proof}
First, we note that part (1) is immediate from the construction of $h$ since we showed that $\Ex[g(Z_1,Z_2)]=\binom{n}{2}/2$ when finding $\Ex[I(\sigma)]$. To see part (2), we claim that $h(Z_1,Z_2)\in\{\binom{n}{2}/2,-\binom{n}{2}/2\}$ a.s. If $U_1<U_2$, then we have
\[h(Z_1,Z_2)=\binom{n}{2}\left(\1(X_{1}>X_{2})+\1(A_{1,2})-1/2\right).\]
The events $\{X_1>X_2\}$ and $A_{1,2}=\{X_1=X_2,\text{ both even}\}$ are disjoint. Hence, when $U_1<U_2$, we have $h(Z_1,Z_2)\in\{\binom{n}{2}/2,-\binom{n}{2}/2\}$ a.s. By symmetry (noting $A_{1,2}=A_{2,1}$), we have $h(Z_1,Z_2)\in\{\binom{n}{2}/2,-\binom{n}{2}/2\}$ a.s. when $U_1>U_2$. Since the events $\{U_1<U_2\}$ and $\{U_1>U_2\}$ are disjoint, and $U_1\neq U_2$ a.s., we conclude that $h(Z_1,Z_2)\in\{\binom{n}{2}/2,-\binom{n}{2}/2\}$ a.s. Noting that $\mathbb{E}[h(Z_1,Z_2)]=0$, we conclude that $\sigma^2:=\Var(h(Z_1,Z_2))=\binom{n}{2}^2/4$, which gives part (2). For part (3), we define
\begin{equation}
h_1(x_1,u_1)=\mathbb{E}[h((X_1,U_1),(X_2,U_2))\mid X_1=x_1,U_1=u_1].
\end{equation}
This allows us to find that 
\begin{equation}
\begin{aligned}
    h_1(x_1,u_1)=\binom{n}{2}\bigg(\mathbb{P}(X_2<x_1)&\mathbb{P}(U_2>u_1)+\mathbb{P}(X_2>x_1)\mathbb{P}(U_2<u_1)\\
    &+\mathbb{P}(U_2>u_1)\mathbb{P}(X_2=x_1,\text{ both even})\\&+\mathbb{P}(U_2<u_1)\mathbb{P}(X_2=x_1,\text{ both even})-\frac{1}{2}\bigg).\\
\end{aligned}
\end{equation}
Computing probabilities, we have
\[h_1(x_1,u_1)=\binom{n}{2}\left(\frac{x_1-1}{2m}(1-u_1)+\frac{2m-x_1}{2m}u_1+\frac{\1(x_1\text{ even})}{2m}-\frac{1}{2}\right).\]
Therefore,
\begin{equation}
\label{eq:h1z1}
    \begin{aligned}
        h_1(Z_1)&=h_1(X_1,U_1)\\
        &=\binom{n}{2}\left(\frac{X_1-1}{2m}(1-U_1)+\frac{2m-X_1}{2m}U_1+\frac{\1(X_1\text{ even})}{2m}-\frac{1}{2}\right).
    \end{aligned}
\end{equation}
To compute $\Var(h_1(Z_1))$, we let
\[Y=\frac{X_1-1}{2m}(1-U_1)+\frac{2m-X_1}{2m}U_1+\frac{\1(X_1\text{ even})}{2m}-\frac{1}{2}.\]
We can rewrite
\[Y=\frac{X_1-1}{2m}+\frac{2m-2X_1+1}{2m}U_1+\frac{\1(X_1\text{ even})}{2m}-\frac{1}{2}.\]
Conditioning on $X_1$, we have 
\[\mathbb{E}[Y\mid X_1]=\frac{X_1-1}{2m}+\frac{2m-2X_1+1}{4m}+\frac{\1(X_1\text{ even})}{2m}-\frac{1}{2}=-\frac{1}{4m}+\frac{\1(X_1\text{ even})}{2m}\]
and
\begin{equation}
    \begin{aligned}
    \Var(Y\mid X_1)&=\left(\frac{2m-2X_1+1}{2m}\right)^2 \cdot\Var(U_1)\\
    &=\frac{(2m-2X_1+1)^2}{48m^2}.
    \end{aligned}
\end{equation}
By the conditional variance formula, we have
\begin{equation}
    \begin{aligned}
        \Var(Y)&=\Var(\mathbb{E}[Y\mid X_1])+\mathbb{E}[\Var(Y\mid X_1)]\\
        &=\Var\left(-\frac{1}{4m}+\frac{\1(X_1\text{ even})}{2m}\right)+\mathbb{E}\left[\frac{(2m-2X_1+1)^2}{48m^2}\right]\\
        &=\frac{1}{16m^2}+\frac{4m^2-1}{144m^2}\\
        &=\frac{m^2+2}{36m^2}
    \end{aligned}
\end{equation}
where the third equality holds because $\1(X_1\text{ even})\sim\text{Bernoulli}(1/2)$ and $X_1$ is discrete uniform on $[2m]$. We conclude that
\[\sigma_1^2:=\Var(h_1(Z_1))=\binom{n}{2}^2\Var(Y)=\binom{n}{2}^2\frac{m^2+2}{36m^2}.\]
This proves (3). Lastly, we show (4). Recalling the expression for $h_1(Z_1)$ in Equation $\eqref{eq:h1z1}$, the triangle inequality implies that
\begin{equation}
    \begin{aligned}
        |h_1(Z_1)|&\leq\binom{n}{2}\left(\left|\frac{X_1-1}{2m}\right||1-U_1|+\left|\frac{2m-X_1}{2m}\right||U_1|+\left|\frac{\1(X_1\text{ even})}{2m}\right|+\frac{1}{2}\right)\\
        &\leq\binom{n}{2}\left(\frac{2m-1}{2m}+\frac{2m-1}{2m}+\frac{1}{2m}+\frac{1}{2}\right)=\binom{n}{2}\left(\frac{5}{2}-\frac{1}{2m}\right)\leq\frac{5}{2}\binom{n}{2}.
    \end{aligned}
\end{equation}
It follows from monotonicity that $\mathbb{E}[|h_1(Z_1)|^3]\leq\left(\frac{5}{2}\right)^3\binom{n}{2}^3$, which proves (4). 
We now have all of the pieces of the puzzle to obtain Theorem $\ref{CLT inversions}$.
\end{proof}
\begin{proof}[Proof of Theorem $\ref{CLT inversions}$]
Let the random variables $Z_1,\ldots, Z_n$ and the function $g$ be as in Proposition $\ref{prop:inversionsmaindecomp}$. Then we have the decomposition
\[I(\sigma)=_d\frac{1}{\binom{n}{2}}\sum_{1\leq i<j\leq n}g(Z_i,Z_j).\]
Letting $h(Z_i,Z_j)=g(Z_i,Z_j)-\binom{n}{2}/{2}$, we have that by construction, $h$ is a symmetric, Borel measurable function and the random variables $Z_1,\ldots, Z_n$ are i.i.d. Hence, the random variable $W_n$ from Equation $\ref{eq:wn}$ is a $U$-statistic. Denoting $Z\sim\mathcal{N}(0,1)$, we obtain from applying Proposition $\ref{prop:U-statistics}$ to the data in Proposition $\ref{prop:dataforclt}$ that
\begin{equation}
    \begin{aligned}
        d_K\left(\frac{\sqrt{n}}{2\binom{n}{2}\sqrt{\frac{m^2+2}{36m^2}}}\left(\underbrace{I(\sigma)-\frac{\binom{n}{2}}{2}}_{W_n}\right),Z\right)&\leq\frac{7\cdot\left(\frac{5}{2}\right)^3\binom{n}{2}^3}{\sqrt{n}\binom{n}{2}^3\left(\frac{m^2+2}{36m^2}\right)^{3/2}}+\frac{(1+\sqrt{2})\frac{\binom{n}{2}}{2}}{(2(n-1))^{1/2}\binom{n}{2}\left(\frac{m^2+2}{36m^2}\right)}\\
        &=\frac{216\cdot7\cdot\left(\frac{5}{2}\right)^3}{\sqrt{n}\left(1+\frac{2}{m^2}\right)^{3/2}}+\frac{9\sqrt{2}\cdot(1+\sqrt{2})}{\sqrt{n-1}\left(1+\frac{2}{m^2}\right)}\\
        &\leq\frac{216\cdot7\cdot\left(\frac{5}{2}\right)^3+18\cdot(1+\sqrt{2})}{\sqrt{n}}
    \end{aligned}
\end{equation}
as desired, since the numerator is a constant which does not depend on $n$ and $m$.
\end{proof}

\section*{AI Disclosure}
The author acknowledges that ChatGPT 5.5 was used to assist with typesetting material and the proof of Lemma $\ref{lem:lonelycoupling}$.
\bibliographystyle{amsplain}
\bibliography{AC}
\end{document}